\documentclass[]{article}
\usepackage{amssymb,amsmath,hyperref,verbatim,cite,graphicx}
\usepackage{float}
\usepackage{authblk}
\newfloat{Chart}{htbp}{grf}
\makeatletter \let\cl@chapter\relax \makeatother
\usepackage{multirow, longtable}
\usepackage{booktabs}
\usepackage{diagbox}
\usepackage{comment}
\usepackage{pdflscape}
\usepackage{mathtools}

\usepackage{longtable}
\usepackage{enumerate}
\usepackage[dvipsnames]{xcolor}
\usepackage{tikz}
\tikzstyle{startstop} = [rectangle, rounded corners, minimum width=3cm, minimum height=1cm,text centered, draw=black, fill=red!30]
\tikzstyle{io} = [rectangle, rounded corners, minimum width=3cm, minimum height=1cm,text centered, draw=black, fill=green!30]
\tikzstyle{process} = [rectangle, minimum width=2cm, minimum height=1cm, text centered, draw=black, fill=blue!40]
\tikzstyle{thread} = [rectangle, minimum width=2cm, minimum height=1cm, text centered, draw=black, fill=gray!40]
\tikzstyle{parallel} = [rectangle, minimum width=3cm, minimum height=1cm, text centered, draw=black, fill=orange!40]
\tikzstyle{arrow} = [thick,->,>=stealth]
\usepackage{pgfplots}
\pgfplotsset{width=10cm,compat=1.9}
\usepackage[ruled,vlined,linesnumbered]{algorithm2e}
\usepackage{hyperref}
\usepackage{subcaption}
\hypersetup{
    colorlinks=true,
    linkcolor=blue,
    filecolor=magenta,      
    urlcolor=cyan,
}
\usepackage{cleveref, xurl}

\begin{document}

\newcommand{\red}{\color{red}}
\newcommand{\blue}{\color{blue}}
\newcommand{\green}{\color{green}}
\allowdisplaybreaks

\title{Serial and Parallel Two-Column Probing for Mixed-Integer Programming}
\author[]{Yongzheng Dai}
\author[]{Chen Chen}
\affil[]{ISE, The Ohio State University, Columbus, OH, USA}

\maketitle

\begin{abstract}
Probing in mixed-integer programming (MIP) is a technique of temporarily fixing variables to discover implications that are useful to branch-and-cut solvers. Such fixing is typically performed one variable at a time---this paper develops instead a two-column probing scheme that instead fixes a pair of variables per iteration. Although the scheme involves more work per iteration compared to the one-column approach, stronger implied bounds as well as more conflicts identified may compensate. Indeed, our prototype implementation was awarded first prize at the MIP Workshop 2024 Computational Competition on novel presolving approaches. This paper presents the aforementioned (serial) prototype and additionally develops an efficient parallelization, leveraging hardware acceleration to further improve overall solve times. Compared to serial two-column probing, our parallel version sacrifices some strength per-pair probed in exchange for greatly increasing the total number of such probings; computational experiments demonstrate its promise.
\end{abstract}

\section{Introduction}
For mixed integer programming (MIP) problems, presolving (or preprocessing) is a set of techniques run prior to the main branch-and-bound search that aims to transform a given instance into an equivalent but (hopefully) easier-to-handle formulation.  Techniques include eliminating redundancies and strengthening the formulation \cite{achterberg2020presolve,PresolveProgress}. Presolve plays an important role in practical results for MIP solvers \cite{achterberg2013mixed}, as well as linear programming (LP) \cite{PresolveLP,PresolveInteriorPoint}, and mixed integer nonlinear programming (MINLP) \cite{PresolveMINLP,PresolveMISDP}. Besides general solvers, there are also problem-specific preprocessing methods (see, e.g. \cite{application1,application2,application3}).

Probing plays an important role in MIP presolves with demonstrated impact on solvers, e.g. \cite{achterberg2020presolve}. This technique temporarily fixes binary variables to obtain a variety of useful implications such as variable bound strengthening, and aggregations \cite{FirstProbing,Achterberg2005,scip, aggregation,hoffman1993solving,TwoRowColumn,achterberg2007conflict,achterberg2020presolve}. The MIP version of probing follows probing for SAT solvers \cite {lynce2003probing}, and has subsequently also found use in extensions such as MINLP \cite{ProbMINLP} and multi-stage integer programming \cite{ProbStochastic, ProbStochastic2}. 

The standard approach to probing involves fixing one binary variable at a time and running e.g. domain propagation \cite[Sec. 7.2]{achterberg2020presolve} and \cite[Sec. 10.6]{scip}. Probing all binary variables can be potentially time-consuming due to the quadratic runtime complexity associated with the number of nonzero constraint coefficients; thus, solvers typically set a workload-limit and attempt to prioritize variables with a probing order \cite{scip,scip1,scip2,achterberg2020presolve,gleixner2023papilo}. In this paper, we propose a new two-column probing approach, involving pairs of variable per iteration. As with standard probing, in two-column probing we rely on workload-limits and judicious selection of (pairs of) variables.  Our approach is a natural generalization of one-column probing that follows a similar philosophy applied to other parts of MIP such as: multiple columns to reduce models \cite{ParallelRowColumn,achterberg2020presolve}; multiple-row reduction \cite{MultipleRow}; matrix sparsification \cite{MatrixSparse1,PresolveInteriorPoint,achterberg2020presolve}; and two-row/column bound strengthening \cite{TwoRowColumn}. We highlight  that two-variable probing was considered in Atamturk et al. \cite{atamturk2000conflict} (see also \cite{atamturk2000mixed}), deployed in the context of conflict graph construction---pairs of variables used to detect infeasible and also suboptimal assignments.  In this paper we apply probing in a modern, integrated fashion: not only to detect conflicts but also can be propagated repeatedly for aggregations, bound strengthening and implications. Indeed, the clique table from a conflict graph procedure can be used to warmstart our probing orders.

Additionally, we develop a parallelization scheme for two-column probing that enables substantially more pairs to be processed, but may generate weaker (local) per-probing results compared to the (global) serial scheme. Parallelization has become increasingly important component of MIP solvers \cite{bixby2007progress,berthold2012solving,achterberg2013mixed,koch2022progress}; we note also a growing literature in massively parallel computing for MIP (e.g. \cite{phillips2006massively,eckstein2015pebbl,shinano2016solving,koch2012could,shinano2018fiberscip,perumalla2021design}). Notably, Gleixner, Gottwald, and Hoen \cite{gleixner2023papilo} parallelized one-column probing by distributing each task of probing a variable to different threads with a similar scheme. We note also that this paper shares the philosophy of our preceding work on parallel conflict graph management \cite{parallelclqmerge}: to develop intensive versions of existing preprocessing methods for which the overhead can be mitigated by parallel computing with carefully managed communication costs.  


The remainder of the paper is organized as shown in Chart~\ref{fig:org_chart}. Sec. \ref{sec:notation} introduces the notation, and Sec. \ref{sec:simplepresolve} presents some simple pre-presolves we used before conducting two-column probing. After simple pre-presolves, there are serial two-column probing shown in Sec.\ref{sec:serialprobing} and the parallel version in Sec.\ref{sec:parallelprobing}. Sec. \ref{sec:experiment} describes numerical experiments. Sec. \ref{sec:conclusion} concludes.

\begin{Chart}[!htbp]
\centering
    \begin{tikzpicture}[scale=.75, node distance = 1.5cm]
        \node (start) [startstop] {Sec. \ref{sec:notation}: Notation};
        \node (simplepresolve) [io, below of=start] {Sec. \ref{sec:simplepresolve}: Simple Pre-Presolves};
        \node (serialprob) [process, below of=simplepresolve, xshift=-3cm] {Sec. \ref{sec:serialprobing}: Serial Two-Column Probing};
        \node (parallelprob) [process, below of=simplepresolve, xshift=3cm] {Sec. \ref{sec:parallelprobing}: Parallel Two-Column Probing};
        \node (experiments) [io, below of=parallelprob, xshift=-3cm] {Sec. \ref{sec:experiment}: Numerical Experiments};
        \node (conclusion) [startstop, below of=experiments] {Sec. \ref{sec:conclusion}: Conclusion};
        \draw [arrow] (start) -- (simplepresolve);
        \draw [arrow] (simplepresolve) -- (serialprob);
        \draw [arrow] (simplepresolve) -- (parallelprob);
        \draw [arrow] (serialprob) -- (experiments);
        \draw [arrow] (parallelprob) -- (experiments);
        \draw [arrow] (experiments) -- (conclusion);
    \end{tikzpicture}
    \caption{Organization of the Paper}
    \label{fig:org_chart}
\end{Chart}

\section{Setup}
\subsection{Notation and Basic Definitions}\label{sec:notation}
Consider a mixed integer program (MIP) in the following generic form:
\begin{equation}\label{eq:MIP}
    \begin{aligned}
        \min\ &c^Tx\\
        \mbox{s.t.}\quad &Ax \circ b\\
        &\ell \leq x\leq u\\
        &x_j\in \mathbb{Z} \mbox{ for all } j \in \mathcal{I}
    \end{aligned}
\end{equation}
with parameters $A\in \mathbb{R}^{m\times n}$, $c\in \mathbb{R}^n$, $b\in \mathbb{R}^m$, $\ell \in (\mathbb{R}\cup {-\infty})^n$, and $u \in (\mathbb{R}\cup {\infty})^n$. Variables $x\in \mathbb{R}^n$, and $x_j\in \mathbb{Z}$ for $j\in \mathcal{I} \subseteq \mathcal{N} = \{1,...,n\}$. Constraints with relations $\circ_i\in\{=,\leq,\geq\}$ for each row $i\in \mathcal{M} = \{1,...,m\}$. Furthermore, let $\mathcal{B} := \{j\in \mathcal{I}\ |\ 0\leq x_j\leq 1\}$ be the set of indices for all binary variables. Let $NNZ$ denote the number of nonzero coefficients in $A$.

A \emph{coupling matrix} $CM\in\mathbb{Z}_+^{|\mathcal{B}|\times |\mathcal{B}|}$, where $CM_{i,j}$ is the number of constraints containing binary variables $x_i$ and $x_j$ with $i, j\in\mathcal{B}$, and $CM_{i,i} = 0$. A set packing (clique) constraint is defined as
\begin{equation*}
    \sum_{j\in \mathcal{B}_S}x_j \leq 1,
\end{equation*}
for some $\mathcal{B}_S\subseteq \mathcal{B}$. Set packing constraints can, in turn, be inferred from the general MIP formulation. For a given mixed-integer constraint $A_{i\cdot}x \circ_i b_i$, if $\circ_i$ is $\leq$, then we can extract a pure binary variables constraint (PBC) denoted as $\mbox{PBC}(i)$:
\begin{equation}\label{eq:pure_binary}
    \sum_{j\in \mathcal{B}: A_{ij}>0}A_{ij}x_j - \sum_{j\in \mathcal{B}: A_{ij}<0}A_{ij}\bar{x}_j \leq b_i - \inf\{\sum_{j\not\in \mathcal{B}}A_{ij}x_j\} - \sum_{j\in \mathcal{B}: A_{ij}<0}A_{ij}.
\end{equation}
Note that if $\circ_i$ is $\geq$, we can rewrite the constraint as $-A_{i\cdot}x \leq -b_i$, and if $\circ_i$ is $=$, we can split the constraint into two constraints $A_{i\cdot}x \leq b_i$ and $-A_{i\cdot}x \leq -b_i$. So we only consider constraints with the form of $A_{i\cdot}x \leq b_i$ in the remainder of this section. Let 
\[\mathcal{C}:=\{i\ |\ i\in \mathcal{M}, \mbox{ PBC}(i) \mbox{ is a clique constraint}\},\]
which called as \emph{clique index set}. Let $\mathcal{S}_{sp}:=\{\mbox{PBC}(i)\ |\ i\in\mathcal{C}\}$, which is a collection of all clique constraints, while $\mathcal{S}_{nsp}:=\{A_{i\cdot}x\circ_i b_i\ |\ i\not\in\mathcal{C}\}$ as a complementary set. A \emph{clique table} $CT$ includes $|\mathcal{B}|$ sets, where 
\[CT[j]:=\{i\in\mathcal{M}\ |\ A[i,j] \neq 0, j\in \mathcal{C}\}.\]
The \emph{conflict number}
\begin{equation}\label{eq:coupnumber}
    conf_{j} = \sum_{i\in CT[j]}\|A[i,:]\|_0,
\end{equation}
is the summation of numbers of nonzero elements in constraints $i\in CT[j]$, where $CT[j]$ represents the cliques containing variable $j$.  The conflict number is used to estimate the number of conflicts associated with $x_j$. The aggregation form $AF$ stores relations like $x_k = ax_j + b$, and the implication graph $IG$ stores relations like $x_j = 1 \rightarrow x_k \geq b$.

\subsection{Simple Pre-Presolves}\label{sec:simplepresolve}
At the preprocessing stage, given a $\mbox{MIP} = (\mathcal{M}, \mathcal{N}, \mathcal{I}, A, b, c, \circ, \ell, u)$, we perform some pre-presolves to detect set packing constraints and constructing $CM$ and $CT$ for subsequent probing, which is described in Alg~\ref{Alg:simple_presolve} and Alg~\ref{Alg:build_matrices}. This applies well-studied techniques (see e.g. Achterberg et al. \cite{achterberg2020presolve}).

\begin{algorithm}[ht]
\SetAlgoLined
\LinesNumbered
\SetKwRepeat{Do}{do}{while}
\SetKwInput{Input}{Input}
\SetKwInput{Output}{Output}
\Input{$\mbox{MIP} = (\mathcal{M}, \mathcal{N}, \mathcal{I}, A, b, c, \circ, \ell, u)$}
\Output{MIP after applying Simple Presolve, $\mathcal{S}_{sp}$,  $\mathcal{S}_{nsp}$, $\mathcal{C}$}
 Set $\mathcal{S}_{sp}$, $\mathcal{S}_{nsp} :=\emptyset$\;
 Remove empty constraints and singletons from $\mbox{MIP}$ and conduct a one-round single-row bound strengthening to $\mbox{MIP}$\;
 Check integer variable bounds and identify binary variables\;
 \For{$i \in \mathcal{M}$}{
    \eIf{$\mbox{PBC}(i)$ is a set packing constraint}{
    Set $\mathcal{S}_{sp}:= \mathcal{S}_{sp}\cup\{\mbox{PBC}(i)\}$\;
    Set $\mathcal{C} := \mathcal{C}\cup\{i\}$\;}
    {
        Set $\mathcal{S}_{nsp}:= \mathcal{S}_{nsp}\cup\{A_{i\cdot}x\leq b_i\}$\;
    }
 }
 \textbf{return} MIP after applying Simple Presolve, $\mathcal{S}_{sp}$, $\mathcal{S}_{nsp}$, $\mathcal{C}$. 
 \caption{One-round Simple Presolve}
 \label{Alg:simple_presolve}
\end{algorithm}

All techniques mentioned in line 2 of Alg.~\ref{Alg:simple_presolve} can be found in \cite[Sec 3.1 and 3.2]{achterberg2020presolve}. Empty constraints and singletons are discarded as they cannot be probed; moreover, the strengthening of variables bounds allows for detection of binary variables that are otherwise encoded as general integers. The complexity of Alg~\ref{Alg:simple_presolve} is $O(NNZ)$, where $NNZ$ is the number of nonzero elements in $A$ \cite{gleixner2023papilo}.

\begin{algorithm}[ht]
\SetAlgoLined
\LinesNumbered
\SetKwRepeat{Do}{do}{while}
\SetKwInput{Input}{Input}
\SetKwInput{Output}{Output}
\Input{$\mbox{MIP} = (\mathcal{M}, \mathcal{N}, \mathcal{I}, A, b, c, \circ, \ell, u)$, $\mathcal{S}_{sp}$, $\mathcal{S}_{nsp}$, $\mathcal{C}$}
\Output{$CM$, $CT$}
 Set $CM\in \mathbb{Z}_{+}^{|\mathcal{B}|\times|\mathcal{B}|}$ as all-zero sparse matrix, $CT := \cup_{j\in\mathcal{B}}\{\emptyset\}$\;
 \For{$s\in \mathcal{S}_{nsp}$}{
    \If{Length of $s \leq \texttt{size\_limit}$ }{
        \For{all pairs of $(x_i,x_j)$ from $s$}{
            $CM_{i,j} := CM_{i,j} + 1$\;
        }
    }
 }
 \For{$i\in 1:|\mathcal{C}|$}{
    Set $s, j$ as the $i$-th element of $\mathcal{S}_{sp}$ and $\mathcal{C}$ separately\;
    \For{$x_k \in s$}{
        Set $CT[k] := CT[k] \cup \{j\}$\;
    }
 }
 \textbf{return} $CM$, $CT$.
 \caption{Building $CM$ and $CT$}
 \label{Alg:build_matrices}
\end{algorithm}

We use a sparse matrix to represent $CM$, which has quadratic memory complexity associated with the length of cliques. To avoid too much memory consumption, we set a \texttt{size\_limit} in line 3 of Alg.~\ref{Alg:build_matrices} and exclude cliques larger than the \texttt{size\_limit}. We also set a total workload in $CM$ constructing, e.g. the number of nonzero terms processed by our method is at most \texttt{work\_limit}. The specific values of these work-limit parameters in our experiment will be provided in Sec.~\ref{sec:parameter}. In line 10, $s$ is the $i$-th clique in $\mathcal{S}_{sp}$ and $j$ is the $i$-th index in $\mathcal{C}$.

Specifically, because current solvers \cite{gurobi,scip,copt,cplex2009v12} have already implemented Alg.~\ref{Alg:simple_presolve} and it is fast, we excluded its runtime in experiments (which is marginally small regardless). However, we include the runtime of Alg.~\ref{Alg:build_matrices} because different solvers may construct \emph{clique table} (or conflict graph \cite{atamturk2000conflict,achterberg2007conflict}) with different data structures, which makes various building runtime and may impact subsequent probing runtime. And $CM$ may not be used in other presolve routines.

\section{Two-column Probing (serial)}\label{sec:serialprobing}
Chart.~\ref{fig:twocolumn_chart} summarizes our two-column probing method. The remainder of this section provides details on each component.

\begin{Chart}[!htbp]
\centering
    \begin{tikzpicture}[scale=.75, node distance = 1.5cm]
        \node (Input) [startstop] {Input MIP};
        \node (prepresolve) [io, below of=Input] {Pre-Presolves, see Sec. \ref{sec:simplepresolve}, Alg.~\ref{Alg:simple_presolve}};
        \node (buildmatrix) [process, below of=prepresolve] {Construct $CM, CT$, see Sec. \ref{sec:simplepresolve}, Alg.~\ref{Alg:build_matrices}};
        \node (probingorder) [process, below of=buildmatrix] {Determine Probing Order, see Sec. \ref{sec:probingorder}, Alg.~\ref{Alg:prob_order}};
        \node (mainbody) [process, below of=probingorder] {Two-Column Probing, see Sec. \ref{sec:2prob_algorithm}, Alg.~\ref{Alg:twoprob}};
        \node (output) [io, below of=mainbody] {Build Reduced Model, see Sec. \ref{sec:output}};
        \node (Output) [startstop, below of = output] {Output Reduced MIP};
        \draw [arrow] (Input) -- (prepresolve);
        \draw [arrow] (prepresolve) -- (buildmatrix);
        \draw [arrow] (buildmatrix) -- (probingorder);
        \draw [arrow] (probingorder) -- (mainbody);
        \draw [arrow] (mainbody) -- (output);
        \draw [arrow] (output) -- (Output);
    \end{tikzpicture}
    \caption{Two-column Probing}
    \label{fig:twocolumn_chart}
\end{Chart}

\subsection{Selecting Variable Pairs}\label{sec:probingorder}
When considering the order of variables to probe, Achterberg \cite{scip} determines a probing order of binary variables based on three factors: the number of constraints containing a binary variable $x_j$, the number of arcs of $x_j$ in the implication graph, and the number of conflicts between $x_j$ and other binary variables. For our implementation we focus on constructing an independent probing subroutine, and so we do not generate an implication graph prior to running probing; in practice the implication graph is generated from other preprocessing procedures and refined during probing. Thus, we define probing order in a simpler way, estimating the number of conflicts of binary variable $x_j$ by \emph{conflict number} defined in Equation~(\ref{eq:coupnumber}) by $CT$ from Alg.~\ref{Alg:build_matrices}. The number of constraints containing $x_j$ is $\|A[:,i]\|_0$, where the zero-norm represents the number of nonzero elements from one vector. Because pairs of variables are considered, besides the number of constraints containing one binary variable, we consider the number of constraints containing two binary variables. This number is calculated via the \emph{coupling matrix} $CM\in\mathbb{Z}_+^{|\mathcal{B}|\times |\mathcal{B}|}$ from Alg.~\ref{Alg:build_matrices}.

Alg.~\ref{Alg:prob_order} determines probing order. Line 1 selects \texttt{Cand\_Number} candidate pairs of variables via binary search. $conflict$ in lines 3-8 acts a penalty to help avoid conflicting variables in the probing order---conflicts reduce the number of useful combinations among the pair of binary variables, i.e. $(0,0), (0,1),$ $(1,0), (1,1)$. Line 10 fetches up to a fixed number of pairs, \texttt{Max\_Probe\_Number}, for probing, and applies partial sorting (see e.g. \cite{partialsort}).

\begin{algorithm}[ht]
\SetAlgoLined
\LinesNumbered
\SetKwRepeat{Do}{do}{while}
\SetKwInput{Input}{Input}
\SetKwInput{Output}{Output}
\Input{$CM$, $CT$}
\Output{Indices of candidate binary variables for the probing method}
 Set $Cand$ as an index set including indices of largest \texttt{Cand\_Number} cells from $CM$\;
 \For{$(i,j)\in Cand$}{
    \eIf{$CT[i]\cup CT[j] \neq \emptyset$}
    {Set $conflict = 10$\;}
    {
        Set $conflict = 1$\;
    }
    Set $score(Cand_{i,j}) = \frac{10 CM_{ij} + \|A[:, i]\|_0 + \|A[:, j]\|_0 + 3conf_i + 3conf_j}{conflict}$\;
 }
 Partially sort $Cand$ to extract the top \texttt{Max\_Probe\_Number} elements with $score(Cand_{i,j}^{(1)})\geq score(Cand_{i,j}^{(2)})\geq...\geq score(Cand_{i,j}^{(1000)})$\;
 \textbf{return} Top $10^3$ elements extracted in Line 5.
 \caption{Sort Variables for Probing}
 \label{Alg:prob_order}
\end{algorithm}

\subsection{Two Column Probing} \label{sec:2prob_algorithm}
After selecting index pairs, we probe every selected pair of binary variables with Alg.~\ref{Alg:twoprob}.

\begin{algorithm}[ht]
\SetAlgoLined
\LinesNumbered
\SetKwRepeat{Do}{do}{while}
\SetKwInput{Input}{Input}
\SetKwInput{Output}{Output}
\Input{MIP, candidate pairs of variables $Cand$, $CT$, $ub, lb$}
\Output{$ub, lb$, $AF$, new conflicts set $New\_Conflict$}
 Set $AF, IG, New\_Conflict:=\emptyset$\;
 \For{$(x_i,x_j)\in Cand$}{
    If both $x_i$ and $x_j$ have been probed, go to next iteration\;
    \For{$v$ = $(0,0),(1,0),(0,1),(1,1)$ }{
        Fix $x_i,x_j = v_1, v_2$\;
        Propagate $CT, AF, IG$ to obtain temporary variable bounds $ub_{temp}, lb_{temp}$\;
        Call domain propagation to obtain implied bounds $ub(v), lb(v)$\;
        If the problem is infeasible, add $(2v_1-1)x_i + (2v_2-1)x_j \leq v_1+v_2-1$ to $CT$ and $New\_Conflict$\;
     }
     If there are two or more values of $v$ are infeasible, fix $x_i$ or $x_j$\;
     Update $ub = \min\{ub, \max\{ub(v)\}\}$, $lb = \min\{lb, \max\{lb(v)\}\}$\;
     Update $AF, IG$ for $x_i$ and $x_j$ separately\;
     Update $IG$ for $(x_i, x_j)$\;
     Check termination condition\;
 }
 \textbf{Return} $ub, lb$, $AF$, $New\_Conflict$.
 \caption{Two-column Probing}
 \label{Alg:twoprob}
\end{algorithm}

In the input, $Cand$ includes \texttt{Max\_Probe\_Number} pairs of binary variables selected from Alg.~\ref{Alg:prob_order}. Line 7 applies domain propagation (see e.g. \cite[Sec. 3.2]{achterberg2020presolve}) to get implied variable bounds. In line 8, $(2v_1-1)x_i + (2v_2-1)x_j \leq v_1 + v_2 -1$ cuts off $x_i = v_1, x_j = v_2$, which is a detected clique. Line 10 fixes variables as follows:
\begin{itemize}
    \item Fixes $x_i = 1$ if both $x_i = 0, x_j = 1$ and $x_i = 0, x_j = 0$ are infeasible. Likewise with the complement, $x_i = 0$.
    \item Fixes a pair of a variables if the other 3 of 4 assignments are infeasible.  For instance,  if $(x_i, x_j) \in {(0,0), (1, 0), (0,1)}$ are infeasible then $x_i = 1, x_j = 1$.
    \item Replace $x_j$ with $x_i$ in the formulation if both $x_i = 0, x_j = 1$ and $x_i = 1, x_j = 0$ are infeasible (i.e.  $x_i = x_j$). Likewise check if this holds for the complement: $x_i = 1 - x_j$. 
\end{itemize}

Line 11 updates global variable bounds. In line 12, we update $AF$ and $IG$ for each variable in a pair separately. For a given variable $x_i$ we set, 
\[ub(x_i = 1) := \max\{ub(1,1), ub(1,0)\}, lb(x_i = 1) := \min\{lb(1,1), lb(1,0)\},\]
\[ub(x_i = 0) := \max\{ub(0,1), ub(0,0)\}, ub(x_i = 0) := \min\{lb(0,1), lb(0,0)\}.\] 
Aggregation and implication are then handled as in the standard one-variable probing \cite[Sec. 7.2]{achterberg2020presolve} or \cite[Alg. 10.12, step 4 (c)(d)]{scip}. After line 12, we also add implications like $x_i = 1, x_j = 1 \rightarrow x_k \geq b$ to $IG$. To avoid redundancy, pair implications are only added if the individual implications $x_i = 1 \rightarrow x_k \geq b$ or $x_j = 1 \rightarrow x_k \geq b$ are not in $IG$. 
 
We defer detailed discussion of the termination condition in line 14 to Sec. \ref{sec:terminate}.

\subsection{Discussion and Analysis} \label{sec:serial_discuss}

\subsubsection{Comparison with One-Column Probing} \label{sec:onecolcomp}  
\textbf{Iteration Complexity. } In terms of computational complexity, within the for-loop in line 2 of  Alg.~\ref{Alg:twoprob}, domain propagation is in $O(NNZ)$, and updating variable bounds, $AF, IG, CT$ are in $O(n)$. Hence, the per-iteration complexity is the same as in one-column probing \cite{gleixner2023papilo}. The maximum number of iterations is $O(|B|^2)$ (governed by the number of possible pairs), whereas standard one-column probing is applied once per variable and so has $O(|B|)$ iterations; in either case, workload limits are imposed to make this a constant factor in practice.  Note that $\lceil \frac{B}{2} \rceil$ pairs in $Cand$ are sufficient to attain the same (or better) results as one-column probing, provided such pairs do not share variables.

Practically speaking, the most time-consuming step in probing (both single- and two-column) is generating $CT, AF, IG$, and domain propagation. Alg.~\ref{Alg:twoprob} calls $CT, AF, IG$, and domain propagation $4$ times, which is twice the computational effort of single-column probing, i.e. running single probing to each $x_i$ and $x_j$ separately. 

\textbf{Strength. } Two-column probing for a pair of variables weakly dominates the one-column equivalent in the following sense. All aggregations and implications from single-column probing will also be detected in line 12. The reason is that the implied bounds in line 12 are no weaker than the implied bounds from the one-column probing. For the variable upper bound, suppose $v_1,v_2\in\{0, 1\}$,
\[ub(x_i = v_1, x_j = v_2) \leq ub(x_i = v_1, x_j\in\{0,1\}),\]
\[ub(x_i = v_1, x_j = 1-v_2) \leq ub(x_i = v_1, x_j\in\{0,1\}),\]
\begin{equation*}
    \begin{aligned}
        \implies &\max\{ub(x_i = v_1, x_j = v_2), ub(x_i = v_1, x_j = 1-v_2)\}\}\\ \leq &ub(x_i = v_1, x_j\in\{0,1\});
    \end{aligned}
\end{equation*}
thus the implied upper bound from two-column-fixing is no weaker than the bound from one-column-fixing. Dominance in the lower bound can be shown similarly. 

The example below shows that two-column probing may return a stronger variable bound in domain propagation.

\textbf{Example} Consider the following system:
\begin{align}
        &t + x_1 + x_2 + x_3 + x_4 \geq 4,\\
        &y-x_1\leq 0,\\
        &y+x_2\geq 1,\\
        &z-x_3\leq 0,\\
        &z+x_4\geq 1,\\
        &y,z \in \{0,1\}, t\in [0, 4], x\in[0,1]^4.
\end{align}
One-column probing probes $y$ and $z$ separately, and the variable bound of $t$ can be strengthened to $[0,3]$. Two-column Probing probes $y$ and $z$ together, and the variable bound of $t$ is updated to $[0,2]$, which is tighter.

Additionally, line 8 in Alg.~\ref{Alg:twoprob} can detect some conflicts between binary variables and line 13 can provide two-variable implications that cannot be obtained by single-column probing. 

\textbf{Balancing Reduction and Runtime.} Although two-column probing can provide possibly stronger (and never weaker) reduction per-pair compared to one-column probing, for practical purposes not all pairs are probed (even in the one-column case only a subset of variables are probed): in line 1 of Alg.~\ref{Alg:prob_order}, only pairs of variables appearing in at least one constraint together are considered. Thus, some binary variables that only appear in constraints alone are ignored by our two-column probing implementation; indeed, in experiments we observe that this constitutes nearly half of all binary variables. Avoiding one-column probes is an implementation choice that allows us to easily test the effects of two-column probing separate from those of one-column probing. The central empirical question, then, is whether one can choose pairs of variables to probe that provide sufficiently strengthened results above-and-beyond the single column (or no probing) counterfactual that can make up for the added compute.

\subsubsection{Termination Conditions} \label{sec:terminate}  
We set three criteria to terminate two-column probing, adapting the one-column procedure from SCIP \cite[Algorithm 10.12]{scip}. First, we set a maximum iteration limit $
\texttt{Max\_Probe\_Number}$, i.e. a limit on the number of pairs probed per thread. Second, we set a global time limit of $30$ seconds. Third, we set a soft stopping rule using a parameter $\texttt{eff} = 0$ to estimate the computational effort expended. At the beginning of each iteration, we set $\texttt{eff} = 0.9\times \texttt{eff}$, and at the end of each iteration, $\texttt{eff} = \texttt{eff} + 110$. If two variables are fixed in one iteration, then $\texttt{eff} = 0$; if only one variable is fixed, set $\texttt{eff} = 0.5\times \texttt{eff}$. If any new conflicts or aggregations are generated, set $\texttt{eff} = 0.8 \times \texttt{eff}$. If we update the $IG$, $\texttt{eff} = 0.9 \times \texttt{eff}$. If the effort score exceeds a threshold \texttt{eff\_threshold}, e.g. $\texttt{eff} > \texttt{eff\_threshold}$, then we stop the algorithm. These parameter values were taken from the analogous dynamic limits in SCIP's one-column probing procedure and were not derived from empirical tuning.

\subsection{Output Management}\label{sec:output}

Two-column probing provides strengthened variable bounds, new conflicts, $AF$, and $IG$. We modify the original MIP by replacing the variable bounds with the strengthened bounds and adding clique cuts from new conflicts and aggregation constraints like $x_i = ax_j +b$ from $AF$. However, the number of implications from $IG$ is quite large, so we discard them; in a more mature implementation, there may be a better management scheme that can pass $IG$ in whole or in part in a productive way. 

\section{Parallelized Two-Column Probing}\label{sec:parallelprobing}
In this section, we parallelize our two-column probing procedure. Our method partitions the binary variables and assigns a partition to each thread, which can then be handled independently (Sec. \ref{sec:varaible_partition}). Locally, two-column probing (Alg.~\ref{Alg:twoprob}) is applied to variables assigned to each thread. Then the results are reduced/merged globally by using a customized implication analysis procedure (Sec. \ref{sec:reduction}). 

A diagram of of our parallel probing is given in Chart.~\ref{fig:paralleltwocolumn_chart}. Note that, if run over the same pairs of variables, the results will tend to be weaker than the serial two-column method: in order to reduce communication costs, our per-thread probing is conducted locally, which could potentially miss out on more effective global updates (i.e. interactions between variables in different threads).  The trade-off between effectiveness per-probing and parallel efficiency (total number of probes conducted) is governed by problem structure.

\begin{Chart}[!htbp]
\centering
    \begin{tikzpicture}[scale=.75, node distance = 1.5cm]
    \path (6,-0.5) node(x) {Determine Probing Order};
    \path (6,-1.25) node(y) {(Sec. \ref{sec:probingorder}, Alg.~\ref{Alg:prob_order}),};
    \path (6, -2) node(z) {Two-Column Probing};
    \path (6, -2.75) node(zz) {(Sec. \ref{sec:2prob_algorithm}, Alg.~\ref{Alg:twoprob}).};
        \node (Input) [startstop] {Input MIP};
        \node (prepresolve) [io, below of=Input] {Pre-Presolves, see Sec. \ref{sec:simplepresolve}, Alg.~\ref{Alg:simple_presolve}};
        \node (buildmatrix) [process, below of=prepresolve] {Construct $CM, CT$, see Sec. \ref{sec:simplepresolve}, Alg.~\ref{Alg:build_matrices}};
        \node (partition) [parallel, below of=buildmatrix] {Partition Variables, see Sec. \ref{sec:varaible_partition}, Alg.~\ref{Alg:var_assign}};
        \node (proborder1) [thread, below of=partition, xshift = -3cm]{Thread 1};
        \node (proborder3) [below of=partition]{$\cdots$};
        \node (proborder5) [thread, below of=partition, xshift = 3cm]{Thread $k$};
        \node (imp_analysis) [parallel, below of=mainbody] {Implication Analysis, see Sec. \ref{sec:reduction}, Alg.~\ref{Alg:imp_analysis}};
        \node (output) [io, below of=imp_analysis] {Build Reduced Model, see Sec. \ref{sec:output}};
        \node (Output) [startstop, below of = output] {Output Reduced MIP};
        \draw [arrow] (Input) -- (prepresolve);
        \draw [arrow] (prepresolve) -- (buildmatrix);
        \draw [arrow] (buildmatrix) -- (partition);
        \draw [arrow] (partition) -- (proborder1);
        \draw [arrow] (partition) -- (proborder3);
        \draw [arrow] (partition) -- (proborder5);
        \draw [arrow] (proborder1) -- (imp_analysis);
        \draw [arrow] (proborder3) -- (imp_analysis);
        \draw [arrow] (proborder5) -- (imp_analysis);
        \draw [arrow] (imp_analysis) -- (output);
        \draw [arrow] (output) -- (Output);
        \draw [black] (6, -3) |- (proborder5);
        \draw[dashed, draw = black, rounded corners=10pt] (3.5, 0) rectangle (8.5, -3);
    \end{tikzpicture}    \caption{Two-column Probing (parallel). Orange color rectangles represent new serial subroutines. White colors are run in parallel, with each thread applying the algorithms in the dashed rectangle.}
    \label{fig:paralleltwocolumn_chart}
\end{Chart}

\subsection{Variable Partitioning} \label{sec:varaible_partition}
The main interdependence between the outer for-loop iterations (line 2 of Alg.~\ref{Alg:twoprob}) of the serial implementation occurs in the inner for-loop (line 4). After setting two binary variables to their implied bounds, two-column probing propagates $CT$ to fix additional binary variables. If these $CT$-fixing variables have been probed in previous iterations, two-column probing may fix more variables via $AF$ and $IG$ associated with these $CT$-fixing variables. Thus, for parallelization, we aim to assign variables that conflict with each other to the same thread. Note, however, that maximizing the number of conflicts assigned to each thread is an NP-hard partitioning problem \cite{chopra1993partition,devine2006partitioning}---hence, we propose a greedy heuristic partitioning method described in Alg.~\ref{Alg:var_assign}.

\begin{algorithm}[!htbp]
\SetAlgoLined
\LinesNumbered
\SetKwRepeat{Do}{do}{while}
\SetKwInput{Input}{Input}
\SetKwInput{Output}{Output}
\Input{$\mathcal{B}$, $\mathcal{S}_{sp}$, thread number $k$}
\Output{Partitioned subsets}
 Set $\mathcal{I}_j := \emptyset$ for $j = 1, ..., k$\;
 Set $assigned_i := 0$ for $i\in \mathcal{B}$\;
 Set $j = 0$\;
 \For{$clq\in \mathcal{S}_{sp}$}{
    \For{$i$-th variables $\in clq$}{
        \If{$assigned_i \leq 0$}{
            Set $\mathcal{I}_{j+1} := \mathcal{I}_{j+1} \cup \{i\}$\;
            Set $assigned_i := 1$, $j := (j + 1)\mod{k}$\;
        }
    }
 }

 \For{$i\in \mathcal{B}$}{
    \If{$assigned_i \leq 0$}{
        Set $j^* := \arg\min\{\|\mathcal{I}_j\|, j = 1,...,k\}$\;
        Set $\mathcal{I}_{j} := \mathcal{I}_{j} \cup \{i\}$\;
    }
 
 }
 \textbf{return} $\mathcal{I}_j$ for $j = 1, ..., k$
 \caption{Assign Variables to Different Threads}
 \label{Alg:var_assign}
\end{algorithm}

Algorithm~\ref{Alg:var_assign} assigns variables from the same clique to the same thread, as shown in Fig.~\ref{fig:var_assign}, where a red dashed rectangle indicates binary variables assigned to one thread. If a variable has already been assigned, the heuristic does not reassign it (lines 6 and 13).

\begin{figure}[!htbp]
\centering
    \begin{tikzpicture}[scale=.75,auto=left]
    \path (-2.5, 5) node(x) {\color{red}Thread 1};
    \path (5.5, -0.5) node(y) {\color{red}Thread 2};
    \path (3.7, 4.3) node(z) {\color{red} Thread 3};
        \node[circle,draw] (a) at (-4,4) {1}; 
        \node[circle,draw] (b) at (-1,4) {2}; 
        \node[circle,draw] (c) at (-4,1) {3};
        \node[circle,draw] (d) at (-1,1) {4};
        \node[circle,draw] (e) at (2,-0.5) {5};
        \node[circle,draw] (f) at (4,1) {6};
        \node[circle,draw] (g) at (4,-2) {7};
        \node[circle,draw] (h) at (1,4) {8};
        \node[circle,draw] (i) at (3.5,2.5) {9};
        \draw[-,black] (a) -- (b);
        \draw[-,black] (a) -- (c);
        \draw[-,black] (a) -- (d);
        \draw[-,black] (b) -- (c);
        \draw[-,black] (b) -- (d);
        \draw[-,black] (c) -- (d);
        \draw[-,black] (d) -- (e);
        \draw[-,black] (e) -- (f);
        \draw[-,black] (e) -- (g);
        \draw[-,black] (f) -- (g);
        \draw[-,black] (h) -- (i);
        \draw[dashed, draw = red, rounded corners=10pt] (-4.5, 4.5) rectangle (-0.5, 0.5);
        \draw[ draw = blue, rounded corners=10pt, rotate=-26] (-2, 1) rectangle (2.5, -0.2);
        \draw[dashed, draw = red, rounded corners=10pt] (1, -0.5) -- (4.5, 2) -- (4.5, -3) -- cycle;
        \draw[dashed, draw = red, rounded corners=10pt, rotate=-31] (-1.8, 4.8) rectangle (2.2, 3.2);
    \end{tikzpicture}
    \caption{Conflict Graph Constructed from 4 Cliques}
    \label{fig:var_assign}
\end{figure}

\subsection{Implication Analysis} \label{sec:reduction}
In one-column probing, a follow-up preprocessing method called implication graph analysis \cite{scip} can use the results of probing to update global bounds and detect new aggregations. Implication analysis is similar to one-column probing but drops domain propagation to save time. However, this procedure is notably less effective for serial two-column probing: based on experiments, we observed that implication graph analysis detects very few new aggregations as our two-column probing tends to probe individual variables multiple times, resulting in diminishing returns. However, our parallel two-column algorithm ignores pairs of variables between two threads, and so we adapt implication analysis by running it only on these neglected pairs to try and recover some of the advantage of global probing. After running local two-column probings, we first update global variable bounds by $ub := \min_{i=1, ..., k}\{ub^i\}$ and $lb := \max_{i=1, ..., k}\{lb^i\}$, where $ub^i, lb^i$ are the upper and lower bounds from the $i$-th thread. We then conduct implication analysis on pairs of variables between threads as described in Alg.~\ref{Alg:imp_analysis}.

\begin{algorithm}[!htbp]
\SetAlgoLined
\LinesNumbered
\SetKwRepeat{Do}{do}{while}
\SetKwInput{Input}{Input}
\SetKwInput{Output}{Output}
\Input{$CT$, $AF$, $IG$, $ub$, $lb$}
\Output{updated $ub$, $lb$, $AF$}
 Set $\mathcal{P}:=\{i\ |\ x_i \mbox{ has been probed}\}$\;
 \For{$i\in \mathcal{P}$}{
    \For{$v = 0, 1$}{
        Propagate $CT$ for $x_i = v$ to fix conflicting variables\;
        Propagate $AF, IG$ for $x_i = v$ and fixed variables to obtain temporary variable bounds $ub_{temp},lb_{temp}$\;
        Update $ub$, $lb$, and $AF$\;
    }
 }
 \textbf{return} $ub$, $lb$, $AF$
 \caption{Implication Analysis}
 \label{Alg:imp_analysis}
\end{algorithm}

For each $x_i, i\in\mathcal{P}$ ($\mathcal{P}$ is defined in line 1 of Alg.~\ref{Alg:imp_analysis}), after setting it to $0$ or $1$, Alg.~\ref{Alg:imp_analysis} follows a procedure similar to Alg.~\ref{Alg:twoprob} but drops domain propagation to save time. In line 4, fixing $x_i$ and propagating $CT$ may fix some variables that were assigned to  different threads, such as nodes $4$ and $5$ in Fig.~\ref{fig:var_assign}. Subsequently, propagating $AF$ and $IG$ uses local probing outcomes across threads to detect potential new/stronger implications for $x_i$. As a result, probing outcomes from different threads are exchanged and merged into the main thread, restoring some of the global effects that were lost by parallelization.
 
\subsection{Discussion and Analysis} \label{sec:parallel_discuss}
Parallel two-column probing is described as Alg.~\ref{Alg:parallel_two_prob} and Chart.~\ref{fig:twocolumn_chart}.

\begin{algorithm}[!htbp]
\SetAlgoLined
\LinesNumbered
\SetKwRepeat{Do}{do}{while}
\SetKwInput{Input}{Input}
\SetKwInput{Output}{Output}
\Input{MIP, thread number $k$}
\Output{Reduced MIP}
 Conduct simple pre-presolve and construct $CT$, $CM$\;
 Assign variables to different threads: $\mathcal{I}_1, ...,\mathcal{I}_k$\;
 \For(parallel){$i\in 1:k$}{
    Determine probing order\;
    Conduct $\texttt{Two-Column Probing}$ for $\{x_j, j\in \mathcal{I}_i\}$\;
 }
 Reduce probing results by $\texttt{Implication Analysis}$ and build reduced MIP\;
 \textbf{return} MIP
 \caption{Parallel Two-Column Probing}
 \label{Alg:parallel_two_prob}
\end{algorithm}

For two-column probing in line 4, we set similar termination conditions discussed in Sec.~\ref{sec:terminate}, but at most probe $\lceil \texttt{Max\_Probe\_Number}/\ln(k+1)\rceil$ pairs of variables, where $k$ is the number of threads. Note that in most cases, the maximum iteration limit and global time limit are relatively loose; thus Alg.~\ref{Alg:parallel_two_prob} can expend similar runtime as the serial version. However, for a given runtime, the parallel version can probe more pairs of variables, as demonstrated by numerical experiments in Sec~\ref{sec:parallel_performance}. As mentioned at the beginning of the section, local probing is (weakly) less effective than global probing for a given pair. Whether the tradeoff in quantity over quality is worthwhile is ultimately an empirical question, which we explore in the next section.

\section{Numerical Experiments}\label{sec:experiment}
This section describes our implementation and experiments.  All code is publicly available at \url{https://github.com/foreverdyz/TwoColumnProbing}.

\subsection{Experimental Setup}

\subsubsection{Test Set}
Experiments are conducted on the benchmark set of the MIPLIB 2017 Collection \cite{miplib}, consisting of $240$ instances. We remove 47 instances for which no two-column probing is available---no pair of binary variables appears in the same constraint---as our subroutine simply terminates with minimal time spent and no effect on the solver. The results below are run on the remaining $193/240$ instances.

\subsubsection{Software and Hardware}
All algorithms are implemented in Julia 1.8.5 \cite{bezanson2017julia}. We solve both original MIPs and reduced MIPs (MIPs after applying the two-column probing) SCIP v0.11.14 \cite{scip1,scip2} and an alpha (prototype) version of JuMP 1.28 \cite{jump} on an Owens node from the Ohio Supercomputer Center. This node includes 2 sockets, with each socket containing 14-core Intel Xeon E5-2680 CPUs. The node has 128 GB of memory and runs 64-bit Red Hat Enterprise Linux 7.9.0. We use $16$ GB, $16$ cores for the parallel computation, and $1$ core for the SCIP solver.

\subsection{Parameter Setting}\label{sec:parameter}
Considering our software and hardware, for Alg.~\ref{Alg:build_matrices}, we set $\texttt{size\_limit} := 10^3$ and $\texttt{work\_limit} := 2\times 10^8$; for Alg.~\ref{Alg:prob_order}, $\texttt{Cand\_Number}:=5\times10^6$ and $\texttt{Max\_Probe\_Number}:=10^3$; the effort score threshold $\texttt{eff\_threshold} : = 1000$ from Sec.~\ref{sec:terminate}.

Instances were predominantly terminated by soft stopping rule.  For example, in Tables.~\ref{tab:org_comp} and~\ref{tab:org_comp_break}, the time limit was hit in either 2/193 or 3/193 instances depending on thread count, with effects on the averages being less than reporting precision  ($0.1\%$). Variable limits were only hit in serial implementation (10/193 instances), with similar negligible impact on analysis. The limited impact is because reaching time and variable limits typically coincided with challenging instances that could not be solved irrespective of probing settings.

\subsubsection{Measurements}

\texttt{Pre Time} measures time spent in two-column probing.  This excludes overhead from the reading time of the input file and the time to build and write the reduced model to file by JuMP, but includes building $CM$, $CT$, selecting variable pairs, two-column probing, and the implication analysis. \texttt{Pairs \#} recorded how many pairs of binary variables are probed by (parallel) two-column probing.

\texttt{Org Time} is the time reported by SCIP for solving the original models directly read form the input file. \texttt{Rd Time} is the time reported by SCIP for solving the reduced models generated by two-column probing. \texttt{Runtime Comp.} is defined as
\begin{equation}\label{eq:runtime_comp}
    \texttt{Runtime Comp.} := \frac{\texttt{Pre Time} + \texttt{Rd Time}}{\texttt{Org Time}}.
\end{equation}

We also record the primal-dual gap based on the primal and dual bounds reported by SCIP. To get a gap within a fixed range, here we use the primal-dual gap defined by Berthold \cite{BERTHOLD2013611}. This gap is between $[0, 100]$: if SCIP finds an optimal solution or certificates the infeasibility, the gap is $0$; if the solver does not find any feasible solution but also cannot certificate the infeasibility, the gap is $100$; otherwise, the gap is calculated as
\begin{equation}\label{eq:gap}
    \texttt{GAP} := \frac{|\texttt{primal\_bound} - \texttt{dual\_bound}|}{\max\{|\texttt{primal\_bound}|, |\texttt{dual\_bound}|\}}\times 100.
\end{equation}

\subsection{Parallel Performance}\label{sec:parallel_performance}
We present both the geometric mean of probing runtimes (denoted as \texttt{Pre Time}) with a shift of 1 as well the average number of pairs of probed variables (denoted as \texttt{Pairs \#}) in Table~\ref{tab:presolve}. This data is further illustrated in Fig.~\ref{fig:pairs_number}, and we fit the \texttt{Pairs \#} in different number of threads with a sublinear function $f(k) := 35k^{\log_2(1.45)}$ where $k$ is the number of threads and $\log_2(1.45) \approx 0.54 < 1$. 16 threads allowed for approximately 4 times more pairs to be considered.

\begin{table}[!ht]
\centering
\caption{\label{tab:presolve}Parallel Performance: \texttt{Pre Time} is geometric mean with a shift of 1, and \texttt{Pairs \#} is the number of pairs of binary variables processed}
\setlength{\tabcolsep}{2.0mm}{
\begin{tabular}{lcccccc}
\toprule[1pt]
Threads \# &1 Thread  &2 Threads &4 Threads &8 Threads &16 Threads \\
\hline
Pre Time &3.2 &3.4 &3.2 &3.2 &3.1\\
Pairs \# &38.0 &50.2 &74.8 &106.6 &154.5\\
\bottomrule[1pt]
\end{tabular}}
\end{table}

\begin{figure}[!htbp]
\centering
\begin{tikzpicture}[scale=.7]
\begin{axis}[
    axis lines = left,
    xlabel={Number of Threads [$k$]},
    ylabel={\texttt{Pairs \#}},
    xmin=1, xmax=17,
    ymin=30, ymax=160,
    xtick={1,2,4,8,16},
    legend pos=north west,
]
    \addplot[
    color=cyan,
    mark=*,
    style= ultra thick,
    ]
    coordinates {
    (1,38)(2, 50.2)(4, 74.8)(8,106.6)(16, 154.5)
    };
    ]
    \addplot[
    dotted,
    color=purple,
    mark=*,
    style= ultra thick,
    ]
    coordinates {
    (1,35)(2, 50.8)(4, 73.6)(8,106.7)(16, 154.7)
    };
    ]
    \legend{Parallel Two-Column Probing, Fitting Function: $35k^{\log_2(1.45)}$}
\end{axis}
\end{tikzpicture}
\caption{\texttt{Pairs \#} vs. different number of threads.}
\label{fig:pairs_number}
\end{figure}

Table~\ref{tab:presolve} and Fig.~\ref{fig:pairs_number} indicate that (parallel) two-column probing takes $3.1\sim 3.4$ seconds (with a shift of 1), while the number of pairs of probed variables increases with a sublinear rate.

\subsection{Default SCIP}\label{sec:default_scip}
We employed two-column probing with varying numbers of threads (1, 2, 4, 8, 16) to preprocess 193 instances and generate reduced models. Subsequently, we solved both the original and reduced models using SCIP with default settings, a single thread, a memory limit of 16 GB, and a time limit of $3600 - \texttt{Pre Time}$ seconds, where for the original model, time limit is $3600$ seconds. The results of SCIP's solving process are summarized in Table~\ref{tab:org_comp}. 

For column labeled ``serial'', we report for serial two-column probing: in the first row, the geometric mean of \texttt{Runtime Comp.} (defined in Equation~(\ref{eq:runtime_comp})) with a shift of 10 for instances that can be solved by SCIP within the 3600-second time limit on either original or reduced model from serial two-column probing; in the second row, the geometric mean of \texttt{GAP} (defined in Equation~(\ref{eq:gap})) with a shift of 1 for instances that cannot be solved by SCIP within the time limit on either original or reduced model from serial two-column probing. Likewise, for $2, 4, 8, 16$ threads, we report shifted geometric means of \texttt{Runtime Comp.} and \texttt{GAP}.

\begin{table}[!ht]
\centering
\caption{\label{tab:org_comp}Two-Column Probing + Org vs. Org; \texttt{Runtime Comp.} for within time-limit cases and \texttt{GAP} for exceeding time-limit cases}
\setlength{\tabcolsep}{1.2mm}{
\begin{tabular}{lcccccc}
\toprule[1pt]
Threads \# &Serial  &2 Threads &4 Threads &8 Threads &16 Threads \\
\hline
\texttt{Runtime Comp.} &97.6\% &100.0\% &98.3\% &97.5\% &95.4\%\\
\texttt{GAP} &88.4\% &94.2\% &96.5\% &93.8\% &98.2\%\\
\bottomrule[1pt]
\end{tabular}}
\end{table}

From Table.~\ref{tab:org_comp}, serial two-column probing attains a $2.4\%$ speedup compared to SCIP with default settings, as well as substantially improving duality gap on unsolved instances. The parallel method has weaker performance on 2 threads, but provides substantial speedups with 16 threads. 

In Table.~\ref{tab:org_comp_break} we provide runtime breakdowns by problem difficulty, with the brackets set using default SCIP runtimes. Two-column probing appears more advantageous on harder problems.

\begin{table}[!ht]
\centering
\caption{\label{tab:org_comp_break}Two-Column Probing + Org vs. Org; \texttt{Runtime Comp.} for different \texttt{Org Time} brackets.}
\setlength{\tabcolsep}{1.2mm}{
\begin{tabular}{lcccccc}
\toprule[1pt]
Brackets (s) &\# of Instances &Serial  &2 Threads &4 Threads &8 Threads &16 Threads \\
\hline
All &193 &97.6\% &100.0\% &98.3\% &97.5\% &95.4\%\\
$\geq 10$ &186 &96.4\% &98.7\% &95.9\% &96.1\% &93.1\%\\
$\geq 100$ &159 &89.9\% &89.3\% &88.3\% &89.3\% &86.3\%\\
$\geq 1000$ &127 &69.2\% &67.2\% &71.8\% &68.5\% &63.3\%\\
\bottomrule[1pt]
\end{tabular}}
\end{table}

\subsection{Two-Column vs One-Column Probing}\label{sec:disableprob}
We also compare two-column probing with one-column probing disabled (setting parameter ``propagating/probing/maxuseles'' as ``false''). Namely experiments from Sec.~\ref{sec:default_scip} were rerun comparing two-column probing + SCIP without one-column probing vs. default SCIP (including one-column probing). Results are summarized in Table~\ref{tab:no_prob}. 

\begin{table}[!ht]
\centering
\caption{\label{tab:no_prob}Two-Column Probing + Org - One-Column Probing vs. Org; \texttt{Runtime Comp.} for within time-limit cases and \texttt{GAP} for exceeding time-limit cases}
\setlength{\tabcolsep}{1.2mm}{
\begin{tabular}{lcccccc}
\toprule[1pt]
Threads \# &Serial  &2 Threads &4 Threads &8 Threads &16 Threads \\
\hline
\texttt{Runtime} &103.1\% &99.4\% &102.3\% &103.0\% &103.8\%\\
\texttt{GAP} &109.1\% &110.6\% &108.1\% &114.8\% &107.3\%\\
\bottomrule[1pt]
\end{tabular}}
\end{table}


Table.~\ref{tab:no_prob} suggests that standalone two-column probing is less effective than one-column probing as implemented in SCIP. Two contributing factors are suggested here. First, two-column probing affects fewer variables, as binary variables that only appear in constraints alone are not probed. We do this to reduce overlapping efforts with one-column probing. Indeed, we have calculated that approximately half of variables are ignored in these MIPLIB instances.  Second, our two-column implementation is not as well-integrated as SCIP's native one-column method. Unlike our code, SCIP's probing accesses $CT, IG$ and a reduced model that is obtained from other presolve subroutines, and thus probing can be more judiciously be deployed (and with substantially less overhead) by avoiding instances that are more likely to be easy. This is supported by solve time breakouts given in Table~\ref{tab:no_prob_break}.


\begin{table}[!ht]
\centering
\caption{\label{tab:no_prob_break}Two-Column Probing + Org - One-Column Probing vs. Org; \texttt{Runtime Comp.} for different \texttt{Org Time} brackets.}
\setlength{\tabcolsep}{1.2mm}{
\begin{tabular}{lcccccc}
\toprule[1pt]
Brackets &Instances \# &1 Thread  &2 Threads &4 Threads &8 Threads &16 Threads \\
\hline
All &193 &103.1\% &99.4\% &102.3\% &103.0\% &103.8\%\\
$\geq 10$ &186 &104.0\% &98.1\% &101.4\% &103.3\% &104.0\%\\
$\geq 100$ &162 &96.5\% &93.6\% &92.5\% &96.4\% &98.2\%\\
$\geq 1000$ &129 &77.0\% &66.2\% &78.2\% &80.9\% &80.5\%\\
\bottomrule[1pt]
\end{tabular}}
\end{table}


\subsection{Random Seeds with 16 Threads}\label{sec:disableprob}
As a robustness check, we ran additional additional experiments using 5 random seeds per instance to help mitigate solver performance variability \cite{lodi2013performance}. For practical purposes we set a lower total solve time limit of 1800 seconds on these runs, and considered the following configurations:


\begin{itemize}
    \item Org: default SCIP (one-column probing enabled).
    \item Org-1: SCIP with one-column probing disabled.
    \item Two Probing: default SCIP with two-column probing.
    \item Two Probing-1: SCIP with two-column probing but one-column probing disabled.
\end{itemize}

Each percentage in Table~\ref{tab:4runs_time} is the runtime (gap for Table~\ref{tab:4runs_gap}) comparison between the method listed in the column header and the method listed in the row header. For example, the $94.3\%$ in position (1,3) of Table~\ref{tab:4runs_time} is calculated by $\texttt{Two Probing Runtime}/\texttt{Org Runtime}$, i.e. SCIP + (16-threaded) two-column probing is $6.7\%$ faster than default SCIP across $193$ instances $\times 5$ runs. This is in line with the $95.4\%$ single seed result reported in Table~\ref{tab:org_comp}, indicating the effectiveness of two-column probing combined with one-column probing. Likewise, the gap comparison of $96.7\%$ in position (1,3) of Table~\ref{tab:4runs_gap} compares favorably to the single seed result of $98.2\%$. Hence, the results of Section~\ref{sec:default_scip} are supported by these additional experiments. Likewise, results here are consistent with Section~\ref{sec:disableprob}; for example, indicating that two-column probing without one-column probing is slightly less effective than vice versa.

Note that, as with previous subsections, runtime comparison only includes instances that can be solved within the time limit by at least one of the compared methods, while the gap comparison includes instances that cannot be solved within the time limit by at least one method.

\begin{table}[!htbp]
\centering
\caption{\label{tab:4runs_time} \texttt{Runtime Comp.} for within time-limit cases between four methods.}
\setlength{\tabcolsep}{1.7mm}{
\begin{tabular}{l|cccc}
\toprule[1pt]
\texttt{Runtime Comp.} &Org  &Org-1 &Two Probing &Two Probing-1\\
\hline
Org &- &107\% &94.3\% &100.6\%\\
Org-1 &93.5\% &- &92.2\% &94.1\%\\
Two Probing &106.1\% &108.4\% &- &103.8\%\\
Two Probing-1 &99.4\% &106.3\% &96.4\% &-\\
\bottomrule[1pt]
\end{tabular}}
\end{table}

\begin{table}[!htbp]
\centering
\caption{\label{tab:4runs_gap} \texttt{Gap Comp.} for exceeding time-limit cases between four methods.}
\setlength{\tabcolsep}{1.7mm}{
\begin{tabular}{l|cccc}
\toprule[1pt]
\texttt{Gap Comp.} &Org  &Org-1 &Two Probing &Two Probing-1\\
\hline
Org &- &116.7\% &96.7\% &115.0\%\\
Org-1 &85.7\% &- &85.0\% &99.1\%\\
Two Probing &103.4\% &117.7\% &- &117.1\%\\
Two Probing-1 &86.9\% &100.9\% &85.4\% &-\\
\bottomrule[1pt]
\end{tabular}}
\end{table}

\subsubsection{Per-Instance Comparison of Two-Column  vs One-Column}
\begin{figure}[h!] 
    \centering 
    \includegraphics[width=0.75\textwidth]{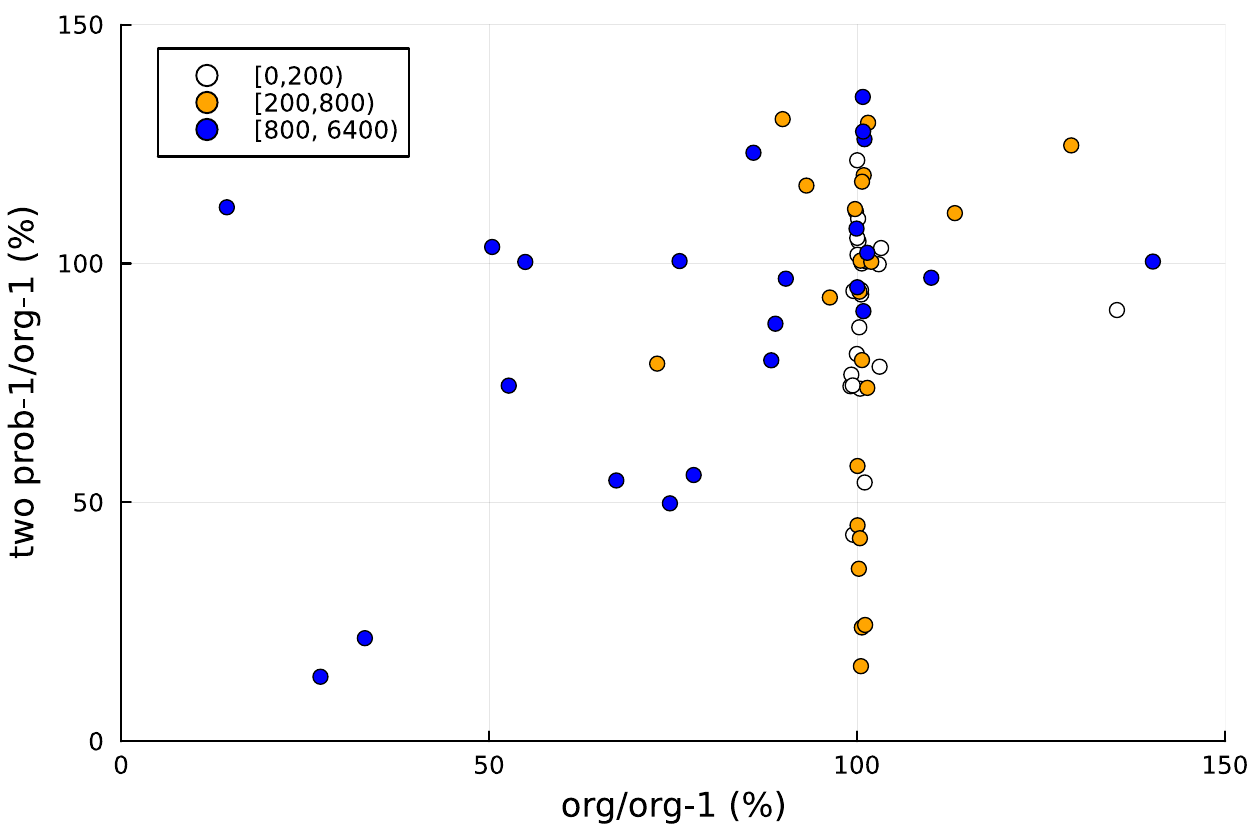}
    \caption{A scatterplot of total runtime effects: vertical axis is two-column probing, horizontal axis is one-column probing. Each point represents five runs of an instance; shade indicates pairs of variables probed.}
    \label{fig:twovsone} 
\end{figure}

Figure~\ref{fig:twovsone} plots instances by runtime comparisons of Two Probing-1/Org-1 on the vertical axis (lower values indicate more time reduction due to two-column probing) and Org/Org-1 on the horizontal axis (lower values indicate more time reduction due to one-column probing). Each point represents a MIPLIB instance (geometric average over 5 random seed runs), with color indicating the number of two-column pairs probed. The figure indicates that when many pairs are probed, the effectiveness of two-column probing on instances is correlated with the effectiveness of one-column probing; otherwise, in instances where we probe fewer pairs, it seems to coincide with minimal one-column probing, hence the cluster of light-coloured points around $100\%$ on the horizontal axis (no one-column effect). This suggests that problem instances that when single-column probing is active, the effects on solve-time would be amplified with the addition of two-column probing. Moreover there are a large number of instances where modest levels of two-column probing is having a substantial effect, but where one-column probing does not.

\subsubsection{Two-Column vs Pairs Probed}
\begin{figure}[h!] 
    \centering 
    \includegraphics[width=0.75\textwidth]{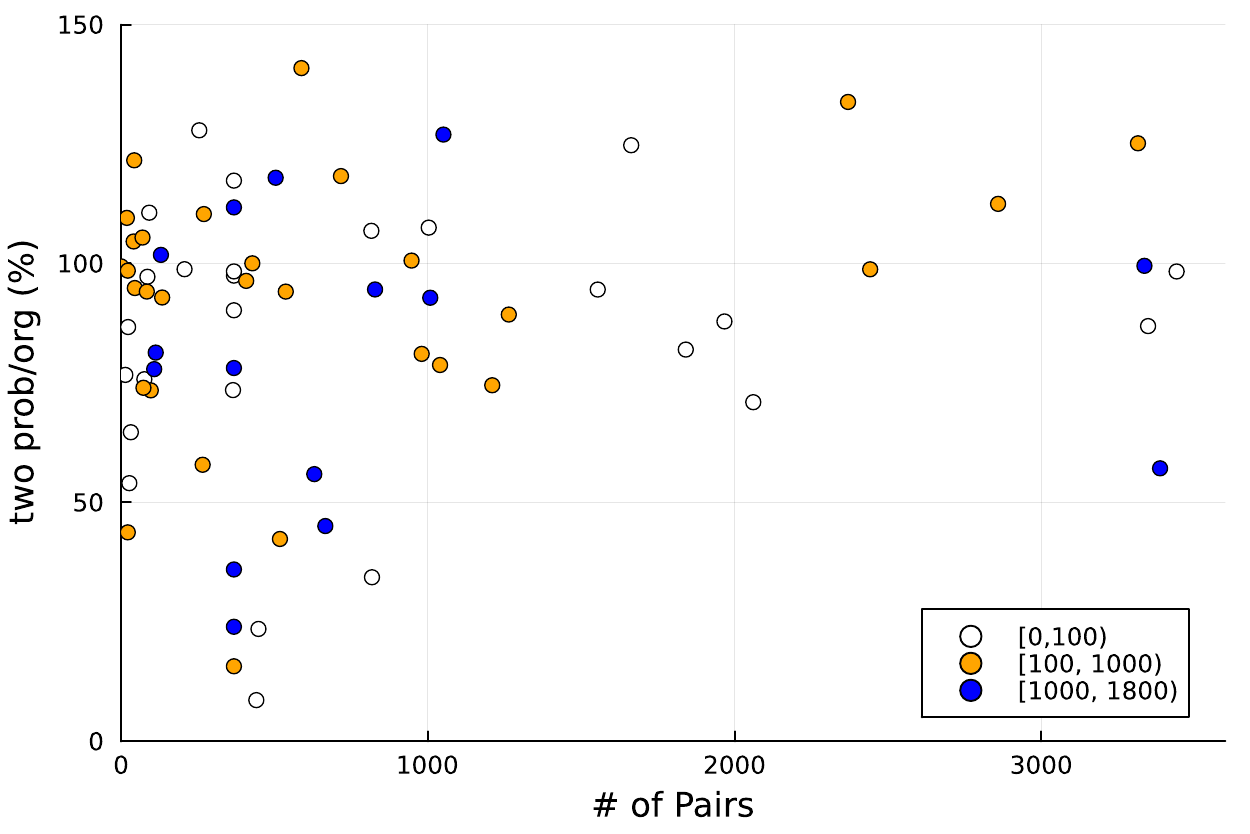}
    \caption{A scatterplot of runtime improvements from two-column probing vs total pairs probed. Each point represents five runs of an instance; shade indicates SCIP total time to solve.}
    \label{fig:twovspairs} 
\end{figure}

Figure~\ref{fig:twovspairs} plots instances by runtime comparisons of Two Probing/Org, i.e. improvement over default SCIP, on the vertical axis (lower values indicate more time reduction from probing) and the number of pairs probed. Each point represents a MIPLIB instance (geometric average over 5 random seed runs), with color indicating time to solve (the longest time to solve among Two Probing and Org configurations). There does not seem to be a clear pattern, although our previous experiments indicate average overall benefits  with higher thread-counts and more pairs probed (see Tables~\ref{tab:org_comp} and ~\ref{tab:org_comp_break}). 

\section{Conclusion}\label{sec:conclusion}
This paper develops a two-column probing procedure, with both a serial and parallel version provided. Although our basic implementation does not have the maturity of fully-integrated one-column probing, empirical results with SCIP nonetheless indicate the promise of our approach, especially on more difficult MIPLIB instances. The parallel version appears advantageous compared with our serial method using 8 or 16 threads. 

For future work, it is worth considering mixed-integer nonlinear problems, where the cost of solving a relaxation tends to be a relatively more expensive compared to the linear setting. In this setting our more intensive version of probing could prove advantageous. 

\section*{Acknowledgements} This work was funded by the Office of Naval Research under grant N00014-23-1-2632.

\bibliographystyle{splncs04}
\bibliography{references}
\end{document}